\documentclass[12pt,oneside,fleqn]{article}

\usepackage[ansinew]{inputenc}
\usepackage[mathscr]{eucal}
\usepackage{amsmath,amssymb,amsthm}
\usepackage{longtable}
\usepackage{color}
\usepackage{cite}
\usepackage{hyperref}
\usepackage{graphicx}
\usepackage{epstopdf}
\usepackage{float}
\usepackage{subcaption}
\usepackage{systeme}
\usepackage{tikz}
\allowdisplaybreaks

\hypersetup{colorlinks, linkcolor=blue, citecolor=blue, urlcolor=blue}

\makeatletter
\long\def\@makecaption#1#2{%
  \vskip\abovecaptionskip\footnotesize
  \sbox\@tempboxa{#1. #2}%
  \ifdim \wd\@tempboxa >\hsize
    #1. #2\par
  \else
    \global \@minipagefalse
    \hb@xt@\hsize{\hfil\box\@tempboxa\hfil}%
  \fi
  \vskip\belowcaptionskip}
\makeatother

\newcommand{\todo}[1][\null]{\ensuremath{\clubsuit}}

\newcommand{\noprint}[1]{}

\newtheorem{theorem}{Theorem}

{\theoremstyle{definition}

\newtheorem{remark}[theorem]{Remark}
\newtheorem*{remark*}{Remark}
}

\newcommand{\checked}[1][\null]{\ensuremath{\boldsymbol{\surd}}}

\newcommand{\R}{\mathbb{R}}

\begin{document}

\par\noindent {\Large\bf
Massively parallel stochastic solution\\ of the geophysical gravity problem
\par}

{\vspace{4mm}\par\noindent {\bf Fabrizio Donzelli$^1$, Alexander Bihlo$^1$, Mauricio Kischinhevsky$^2$ and\\ Colin G.\ Farquharson$^3$
	} \par\vspace{2mm}\par}
{\vspace{2mm}\par\noindent {\it
		$^{1}$~Department of Mathematics and Statistics, Memorial University of Newfoundland,\\ St.\ John's (NL) A1C 5S7, Canada
	}}
	{\vspace{2mm}\par\noindent {\it
			$^{2}$~Institute for Computing, Universidade Federal Fluminense,\\
			Niter\'oi (RJ) 24210-346, Brazil
		}}
	{\vspace{2mm}\par\noindent {\it
			$^{2}$~Department of Earth Sciences, Memorial University of Newfoundland,\\ St.\ John's (NL) A1B 3X5, Canada
		}}
{\vspace{2mm}\par\noindent {\it
				\textup{E-mail:} fdonzelli@mun.ca, abihlo@mun.ca, kisch@ic.uff.br, cgfarquh@mun.ca
}\par}

\vspace{4mm}\par\noindent\hspace*{8mm}\parbox{140mm}{\small
In this paper, we report the advantages of using a stochastic algorithm in the context of mineral exploration based on gravity measurements. This approach has the advantage over deterministic methods in that it allows one to find the solution of the Poisson equation in specified, isolated points without the need of meshing the computational domain and solving the Poisson equation over the entire domain. Moreover, the stochastic approach is embarrassingly parallelizable and therefore suitable for an implementation on multi-core compute clusters with or without GPUs. Benchmark tests are carried out that show that the stochastic approach can yield accurate results for both the gravitational potential and the gravitational acceleration and could hence provide an alternative to existing deterministic methods used in mineral exploration.
}\par\vspace{2mm}

\section{Introduction}\label{sec:intro}

The presence of mineral or oil deposits in the Earth's sub-surface can be detected my measuring small variations of the vertical component of the gravitational acceleration, which are due to changes in the local mass density of the deposit. When geophysicists conduct field studies, they take an average of 20 to 30 gravitational measurements in different locations: at the Earth's surface level, above the surface, or underground. The measurements are then analyzed in order to reconstruct the geometrical and physical properties of the deposit, guided as well by the geological knowledge of the region of study \cite{murray2001best}. Such reconstruction, known as the inverse problem, requires first to have a systematic way to compute the solution of the partial differential equation for the Newtonian potential, the Poisson equation, so as to reconstruct the distribution of the density in the sub-surface.
In this mathematically oriented paper will present a numerical method for solving the direct problem, which involves computing the solution of the Poisson equation for a prescribed sub-surface density structure.  
The geophysical literature proposes numerical solutions of the Poisson equation for gravity which are based on deterministic methods (such as finite differences, finite elements, finite volumes, see e.g.\ \cite{jaha13a,farquharson2009three}).  In this paper we propose instead a probabilistic method, based on the stochastic representation of the solution of the Dirichlet problem associated with the Poisson equation, which is provided by the Feynman--Kac formula \cite{kara91a,aceb05a}.
This formula was introduced by Feynman as an alternative formulation of quantum mechanics, and it was generalized and formulated rigorously by Kac to study parabolic partial differential equations.
Dynkin found a formulation of the formula in the context of elliptic partial differential equations. The Feynman--Kac formula presents the solution of a partial differential equation as an expected value of a certain functional over a collection of random walks, as we will explain in Section \ref{sec:BackgroundMC}.
Despite its theoretical complexity, it admits rather simple numerical approximations as an average over a finite collection of piecewise-linear numerically simulated random paths. We develop and implement a parallel algorithm to construct the finite sample space of paths by means of Monte Carlo simulations.

This paper is meant to be the first in a sequence of papers, in which we will test the accuracy and speed of Monte Carlo based algorithms for solving problems related to exploration geophysics. In particular, we will highlight the following convenient properties of the new approach:

\begin{enumerate}\itemsep=0ex

\item Monte Carlo algorithms do not rely on mesh generation and linear solvers, and therefore are much less memory intensive than deterministic methods.  

\item The performance of the Monte Carlo method is proportional to the number of points at which we are interested in exploring the solution. In practice, gravitational measurements are executed on a small collection of locations. Hence the Monte Carlo algorithm is suitable to model the results of such experiments, since we can adapt it to compute the solution only at the points of interests. 
With deterministic methods, instead, the size of the mesh is going to be rather independent of the number of points of interest, and the performance cannot be improved. 

\item Multiple random walks can be carried out simultaneously using parallel algorithms, therefore speeding up the execution of the simulation, potentially by several orders of magnitude if suitable massively parallel computing infrastructures are available. 

\item The parallelization of the algorithm can be extended by assigning different subsets of points to different parallel processors (or groups of parallel processors), therefore further improving the performance.

\item The low memory requirements and the adaptability of parallelization make Monte Carlo methods highly suitable to be implemented on advanced parallel computing architectures that use GPUs.

\end{enumerate}

In this paper we implement the Monte Carlo method for a rather elementary problem, namely the potential generated by a cube of constant density. It is well known that such a problem admits an analytical solution \cite{waldvogel1976newtonian}: hence, by comparing the results of our numerical method with the analytical solution we have a simple way of testing the accuracy of the stochastic method. In a subsequent work we will propose parallel Monte Carlo simulations for solving the two-dimensional Maxwell's equations modeling magneto-telluric explorations (continuing the work carried out in~\cite{bihl16b}), for which analytical solutions are not generally available. 

The further organization of this paper is as follows. In Section~\ref{sec:BackgroundMC}, we present the necessary mathematical background underlying the stochastic representation of the Poisson equation as relevant for the gravity problem. In Section~\ref{sec:NumericalExperimentsMC}, we present the results of two different implementations of the Monte Carlo algorithm, namely the CUDA-based implementation for the GPU and the MPI-based implementation for a multi-core machine.
In Section~\ref{sec:NumericalExperimentsMC1}, we study the performance of a cluster machine running
the code in parallel across 8 GPUs. Here we also compare accuracy and performance of our Monte Carlo algorithm with a finite element method.  Finally, in Section~\ref{sec:ConclusionsMC} we draw some conclusive remarks and propose future research directions. 

\section{Background}\label{sec:BackgroundMC}

\subsection{The Feynman--Kac formula}

In this section, we recall a few theoretical facts on the Feynman--Kac formula and its numerical implementation through Monte Carlo simulations. 

Let $\Omega\subset\R^3$ be a smooth bounded domain, whose boundary will be denoted as usual by 
$\partial \Omega$. Let $g$ be a continuous function on $\Omega$ and
$f$ be a continuous function on $\partial\Omega$. We are interested in the solution of the following Poisson equation with Dirichlet boundary conditions:
\begin{equation}\label{eq:PoissonEquation}
\Delta u=-g\qquad u_{|_{\partial\Omega}}=f.
\end{equation}
Standard results from functional analysis can be used to prove that this problem is well-posed, which for our practical interest signifies that there exists a unique solution to the above equation.
The key point of the present paper is that the solution of~\eqref{eq:PoissonEquation} has a \textit{stochastic representation}, see e.g.~\cite{kara91a}. It is this stochastic representation that we will use, and 
the construction is as follows. 

Let $\mathbf{x_0}=(x_0,y_0,z_0)\in\Omega$ be the point at which we are interested in computing the solution of the Poisson equation. Let $\mathbf{W}$ be the standard three-dimensional Brownian motion starting at $\mathbf{x}_0$. Intuitively, for $t\geqslant 0$, $\mathbf{W}(t)=(x(t),y(t),z(t))=\mathbf{x}(t)$ is a vector whose 3 components are independent, normally distributed variables, namely $\mathbf{W}\sim \mathcal N(\mathbf{0},I_3)$, with $\mathcal N(\mathbf{0},I_3)$ denoting the distribution of a random $3$-vector of independent standard normal variables, and $\mathbf{W}(0)=\mathbf{x_0}$. In practical terms, each time we simulate a new vector $\mathbf{W}$ we will obtain a different outcome, since $\mathbf{W}$ is a random variable rather than a deterministic vector.  
Let $\boldsymbol{\gamma} (t)$ be the random walk which is the  solution of the stochastic differential equation (practically, $\boldsymbol{\gamma}$ is a ``rescaled" version of $\mathbf{W}$):
\begin{equation}\label{eq:StochasticDE}
\mathrm{d}\boldsymbol{\gamma}(t)=\sqrt{2}\mathrm{d}\mathbf{W}.
\end{equation}

Next, define the first exit time $\tau$ of a random walk $\boldsymbol{\gamma}$ from the domain~$\Omega$ to be the random variable $\tau=\inf_{t>0}\{ t\, |\, \boldsymbol{\gamma} (t)\in\partial\Omega\}$, which represents the first time that the random walk reaches the boundary. 
The Feynman--Kac formula presents the solution of~\eqref{eq:PoissonEquation} at $\mathbf{x}_0$ as:
\begin{align}\label{eq:FeynmanKacFormula}
u(\mathbf{x}_0)=\mathrm{E}_{\mathbf{x}_0}\left(f(\boldsymbol{\gamma}(\tau))+\int_{0}^{\tau(\boldsymbol{\gamma})}g(\boldsymbol{\gamma}(t))\mathrm{d}t\right),
\end{align}
where the expected value~$\mathrm{E}_{\mathbf{x}_0}(\cdot)$ is taken over the sample space of all standard Brownian motions $\mathbf{W}$ in $\Omega$ starting at $\mathbf{x}_0$ (see \cite{kara91a}, Problem 2.25, page 253).

As laid out in measure theory, the expected value of a random variable is computed as an integral over a properly defined measure space. However, the Feynman--Kac formula \eqref{eq:FeynmanKacFormula} is based on the computation of  an integral over a space of paths, which cannot be modeled by a finite-dimensional measure space (see \cite{kara91a}, Chapter 2). 
Even if the above formula is theoretically complex in nature, it offers a method of  evaluating~\eqref{eq:PoissonEquation} numerically which is an alternative to the standard deterministic method (such as a finite difference, finite element or finite volume solver), since we can also evaluate the stochastic representation of the exact solution~\eqref{eq:FeynmanKacFormula} of~\eqref{eq:PoissonEquation} numerically. 
We present the numerical algorithm for \eqref{eq:FeynmanKacFormula} only for 
our specific case of gravity. In this case the partial differential equations is
\begin{align}\label{eq:PoissonEquationNewtonian}
\Delta u=-4\pi G\rho, \quad u(\infty)=0,
\end{align}
where $G$ is Newton's gravitational constant and $\rho$ is the sub-surface density~\cite{jaha13a}. 
To obtain a physically meaningful solution we impose that $u$ vanishes at infinity
(see Section \ref{subsec:boundary} for clarification).

The Feynman--Kac formula in the context of gravity that we will implement below is therefore
\begin{align}\label{eq:FeynmanKacFormulaGravity}
u(\mathbf{x}_0)=4\pi G\,\mathrm{E}_{\mathbf{x}_0}\left(\int_{0}^{\tau(\boldsymbol{\gamma})}\rho(\boldsymbol{\gamma}(t))\mathrm{d}t\right).
\end{align}

The relationship between the gravitational potential~$u$ and the actual gravitational acceleration~$\mathbf{g}$ is given through
\begin{equation}\label{eq:GravitationalAcceleration}
\mathbf{g}=\nabla u.
\end{equation}

Two remarks are necessary here.
Firstly, as was indicated above, the physical potential is obtained by imposing that $u$ vanishes at infinity,
while the solution of equation \eqref{eq:PoissonEquation} is provided for a bounded domain.
In order to overcome this difference we will let $\Omega $ be a ball of sufficiently large radius, whose boundary is sufficiently far from any gravitational source. Vanishing Dirichlet boundary conditions on the solution will therefore be imposed on $\partial \Omega$. An estimate of the systematic error arising from this choice is presented in Section~\ref{subsec:boundary}.

Secondly, another possible source of error arises from the fact that the source term considered in the experiments below is a piecewise constant functions, while the Feynman--Kac formula is known to be applicable, strictly speaking, only to solutions which admit continuous second derivatives. 
This phenomenon will be analyzed numerically in Section~\ref{sec:NumericalExperimentsMC1}.

\subsection{Error due to boundary conditions}\label{subsec:boundary}

The Newtonian potential is the solution of the Poisson equation, with vanishing conditions at infinity. In case of potentials modeled by continuous functions one then assumes that, in Equation \eqref{eq:PoissonEquation}, $\Omega =\R^3$ and $f=0$.
In finding a numerical solution $u_{\rm num}$ one needs necessarily to replace the 3-dimensional space with a bounded but sufficiently large domain $\Omega$, and apply a numerical scheme to solve the equation
\begin{equation}\label{eq:ZeroBoundary}
\Delta u=-4\pi G\rho\qquad u_{|_{\partial\Omega}}=0.
\end{equation}

Hence, any numerical scheme which solves \eqref{eq:ZeroBoundary} will not converge to the solution that represents the true Newtonian potential. In fact, if $\partial \Omega$ is not an equipotential surface for the exact solution, the latter does not vanish identically on the boundary. We can give a rough estimate of the systematic error arising from a non-correct but necessary choice of boundary conditions in the following way. 

Let $u_0$ be the solution of \eqref{eq:ZeroBoundary}, and $u_1$ the Newtonian potential which is the solution of \eqref{eq:PoissonEquationNewtonian}. The difference $u_0-u_1:=u_{\rm e}$ is then the unique solution of 
\begin{equation}\label{equation_error}
\Delta u=0\qquad u_{|_{\partial\Omega}}={u_1}_{|_{\partial\Omega}}.
\end{equation}

Let $B_R$ be the ball of radius $R$ centered at the origin, and denote its boundary sphere by $S_R$.  We recall that the Poisson kernel for the 3-dimensional ball $B_R$  is given by 
\begin{equation}\label{PK}
P_R(\mathbf{x},\boldsymbol{\xi})=\frac{R^2-||\mathbf{x}||^2}{4\pi R||\mathbf{x}-\boldsymbol{\xi}||^3}
\end{equation}
and the solution of \eqref{equation_error} is given by the Poisson integral formula:
\begin{equation}\label{PIF}
u_{\rm e}(\mathbf{x})=\int\limits_{S_R}P(\mathbf{x},\boldsymbol{\xi})u_1(\boldsymbol{\xi})\mathrm{d}S_{\boldsymbol{\xi}}.
\end{equation}

Since in the model problem we are computing values of $u_{\rm e}$ on points very far from $S_R$,
we can approximate $P(\mathbf{x},\boldsymbol{\xi})$ to be a constant, namely 
$P_R(\mathbf{x},\boldsymbol{\xi})\approx \frac{1}{4\pi R^2}$. Moreover, at the boundary sphere, we can approximate the gravitational source by a point mass distribution, hence $u_1(\boldsymbol{\xi})\approx -GM/R$.
The error due to the wrong choice of boundary conditions can therefore be estimated roughly as
$|u_{\rm e}(\mathbf{x})|\approx GM/R$. 
In the experiments we will describe below, this estimate gives as a systematic error of around 1\%, when the computation of the potential
was performed at points which are also sufficiently far from the mass distribution. Other sources of error will dominate instead for
points that are close to the mass distribution, as discussed in Experiment 4, Section~\ref{sec:NumericalExperimentsMC1}.

It is interesting to note however that if the error in the computation of the potential due to the modified boundary conditions is almost constant, it will not affect in a significant way the computation of the acceleration (since the derivative of a constant is equal to zero). 

\subsection{The algorithm}

We fix a point $\mathbf{x}_0\in \Omega$ where we want to compute the numerical solution $u_{\rm num}(\mathbf{x}_0)$ of the Poisson equation~\eqref{eq:ZeroBoundary}. To evaluate~\eqref{eq:FeynmanKacFormulaGravity} numerically, we need to approximate the expected value, the Brownian motion~$\boldsymbol{\gamma}(t)$, and the integral over the sub-surface density distribution. 

The expected value is replaced by a finite sum over~$N$ distinct random walks with associated paths~$\boldsymbol{\gamma}_i(t)$, \
\begin{align*}
{\rm E}_{\mathbf{x}_0}\left( 4\pi G\int_{0}^{\tau(\boldsymbol{\gamma})}\rho(\boldsymbol{\gamma}(t))\mathrm{d}t\right)\approx
\frac{1}{N}\sum_{i=1}^N \left( 4\pi G\int_{0}^{\tau(\gamma_i)}\rho(\boldsymbol{\gamma}_i(t)){\mathrm d}t\right),
\end{align*}
which is the core of the Monte Carlo method.

Next, each stochastic process~$\boldsymbol{\gamma}(t)$, which is governed by the stochastic differential equation
(\ref{eq:StochasticDE}),
where $\mathbf{W}$ is three-dimensional Brownian motion, is approximated by solving a discretized version of the equation. A variety of schemes exists for this purpose such as Taylor-type schemes, exponential time-stepping, and Runge--Kutta methods~\cite{mils04a}. For the sake of simplicity, we choose the simple Euler--Maruyama method here, which discretizes~\eqref{eq:StochasticDE} as
\[
 \boldsymbol{\gamma}^{n+1}=\boldsymbol{\gamma}^n+\sqrt{2\Delta t}\mathbf{W},\qquad \mathbf{W}\sim\mathcal{N}(\mathbf{0},I_3),
\]
where $\Delta t=t^{n+1}-t^n$ is the time step (assumed to be constant throughout the integration), and $\mathcal N(\mathbf{0},I_3)$ is the distribution of random $3$-vectors with independent standard normal variables (i.e.\ with mean zero and unit variance).

Finally, the integral
\[
 I=\int_{0}^{\tau(\boldsymbol{\gamma}_i)}\rho(\boldsymbol{\gamma}_i(t)){\mathrm d}t
\]
in~\eqref{eq:FeynmanKacFormulaGravity} is approximated by a quadrature rule. Here we choose the trapezoidal rule.

With these approximations, each random walk and its corresponding numerical solution are simulated as follows:

\begin{enumerate}
\item Start with $n=0$, $\boldsymbol{\gamma}_i^{n}=x_0$, $I=0$.
\item Simulate a random vector $\mathbf{W}\sim\mathcal{N}(\mathbf{0},I_3)$ and update the position of the process $\boldsymbol{\gamma}_i^{n+1}\leftarrow\boldsymbol{\gamma}_i^{n}+\sqrt{2\Delta t}\mathbf{W}$.
\item Update the Riemann sum, $I\leftarrow I+\frac12(\rho(\boldsymbol{\gamma}_i)^n+\rho(\boldsymbol{\gamma})^{n+1})\Delta t$.
\item Check if $\boldsymbol{\gamma}_i\in\Omega$. If $\boldsymbol{\gamma}_i\in\Omega$, start again at step $2$. Otherwise, record the exit time~$\tau(\boldsymbol{\gamma}_i)$ and compute the solution $u_i=4\pi G I$.
\end{enumerate}

We repeat this procedure~$N$ times and then average over the solutions~$u_i$ obtained for each random walk:
\begin{align}\label{eq:Average}
u_{\rm num}(x_0)=\frac{1}{N}\sum_{i=1}^Nu_i.
\end{align}

While the above algorithm is straightforward, a few remarks are in order to guarantee its effective implementation.

\begin{remark}
A crucial aspect guaranteeing the accuracy of the above algorithm is the correct estimation of the first exit time~$\tau$ of the random process~$\boldsymbol{\gamma}_i(t)$. It is possible that the process started in~$\Omega$ at $t^n$ and finished in~$\Omega$ at $t^{n+1}$ but left the domain in between. To account for this possibility and to estimate the conditional probability of the process exiting during a time step, an interpolating process can be defined, a so-called \textit{Brownian bridge}~\cite{gobe00a}. We have found in practice that the explicit inclusion of such an interpolating process is particularly critical for points~$\mathbf{x}_0$ close to the boundary unless very small time steps~$\Delta t$ are chosen. For points far away from the boundary, we found the difference between using the Brownian bridge or not using it to be negligible.  
\end{remark}

\begin{remark}
As was pointed out in~\cite{aceb05a}, there are three sources of errors in computing the solution to~\eqref{eq:PoissonEquation} via Monte Carlo methods:
(i) replacing the expected value in~\eqref{eq:FeynmanKacFormula} with a finite sum over a finite sample space of random paths~$\boldsymbol{\gamma}_i$, (ii)
the discretization of paths (via discretizing  the stochastic differential equation~\eqref{eq:StochasticDE}) and the replacement of the Riemann integral by a quadrature rule, and (iii) the error in the estimation of the exit time. All these errors are controlled by choosing a sufficiently large number~$N$ of random walks and a sufficiently small time step~$\Delta t$. The precise number of Monte Carlo simulations and required time step size are to be determined based on the accuracy requirements for the resulting numerical solution, which ultimately depends on the accuracy of available field measurements for the gravitational acceleration~$\mathbf{g}$.
\end{remark}

The advantage of using Monte Carlo methods comes from observing that the paths over which we integrate using the Feynman--Kac formula are independent of each other. In practical simulations, this implies that we can carry out separate collection of runs on separate processors or GPU cores, leading to a perfectly scalable algorithm that is suitable for a massively parallel computing infrastructure. This optimal parallel nature of the stochastic algorithm is not shared by standard deterministic methods for solving~\eqref{eq:PoissonEquation}, and is the reason for the renaissance of Monte Carlo methods for the solution of partial differential equations~\cite{aceb05a,aceb07a,aceb10a,bihl14a,bihl15b,bihl16b}, which have long been considered as computationally ineffective due to the exceedingly slow convergence of the Monte Carlo error, which is proportional to~$N^{-1/2}$ when using pseudo-random numbers~\cite{pres07Ay}. Within the framework of parallel computing, the stochastic solution to~\eqref{eq:PoissonEquation} always beats a deterministic solution provided that suitably many compute cores are available.

Once the potential~$u_{\rm num}$ has been computed in the point~$\mathbf{x}_0$ we need to compute the actual gravitational acceleration~$\mathbf{g}$ by numerically solving~\eqref{eq:GravitationalAcceleration}. This is done by computing the solution~$u_{\rm num}$ at points in the neighborhood of~$\mathbf{x}_0$ such that we can compute~$\mathbf{g}$ from a finite difference approximation to~\eqref{eq:GravitationalAcceleration}.

Below we report the results for various experiments using the stochastic algorithm on different machines and with different numbers of processors.

\section{Numerical experiments: GPU vs.\ multi-CPU}\label{sec:NumericalExperimentsMC}

The first two experiments which are presented in this section, were designed to compare the performance of a multi-processor compute cluster with the performance of a single GPU card.  

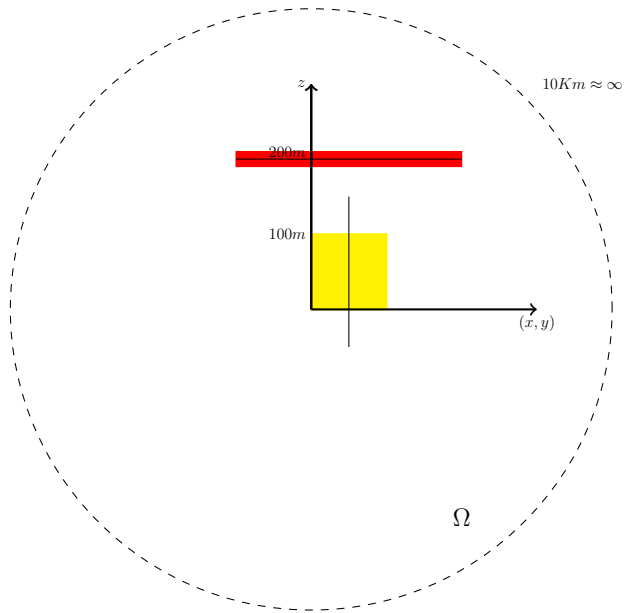
\begin{figure}[!ht]
	\centering
	\begin{tikzpicture}
	\draw[dashed] (0,0) circle (4cm); 
	\filldraw[thick,yellow] (0,0)--(1,0)--(1,1)--(0,1)--cycle; 
	\filldraw[thin,red](-1,2.1)--(-1,1.9)--(2,1.9)--(2,2.1)--cycle; 
	\draw  (-1,2)--(2,2);
	\draw  (0.5,1.5)--(0.5,-0.5);
	\draw[thick,->] (0,0)--(0,3);
	\draw[thick, ->] (0,0)--(3,0) node[anchor=north,scale=0.5]{$(x,y)$};
	\draw[thick, ->] (0,0)--(0,3) node[anchor=east,scale=0.5]{$z$};
	\draw (0,1) node[anchor=east, scale=0.5]{$100m$};
	\draw (0,2.1) node[anchor=east, scale=0.5]{$200m$};	
	\draw (3,3) node[anchor=west,scale=0.5]{$10Km\approx\infty$};
	\draw (2,-3) node[anchor=south,scale=0.8]{$\Omega$};
	\end{tikzpicture}
	\caption{Experiment 1: computation of the potential on the red region. 
		Experiment 2: computation of the potential on along the vertical line.
		The yellow cube represents the mineral deposit.}
\end{figure}\label{f1}

Let us start by describing our simple yet physically relevant test model. The domain $\Omega$ of the Poisson equation is a sphere of radius 10 km. The gravitational acceleration  is generated by a cubical mass $C$ of constant density of 2000 kg/m$^3$, and side length 100 m. In Cartesian coordinates, 
$\Omega$ is the sphere of radius $R=10km$ and centered at the origin, and the cube is the subregion where the coordinates $x,y,z$ range from $0m$ to $100m$.
 We assume that the background potential is zero, 
hence the density $\rho$ from Equation \eqref{eq:PoissonEquationNewtonian} has value $2000kg/m^3$ in $C$, while it vanishes everywhere else.
 We note that in our simulation the cube is positioned in a slightly asymmetric way, with respect to the circular domain~$\Omega$.

As we observed in the previous section, the domain $\Omega$ must be of sufficiently large radius, in order to obtain a realistic approximation of the vanishing condition at infinity that characterizes the Newtonian potential. Such a choice necessarily affects the speed of execution of the Monte Carlo simulation: it is proven in fact that  the expected value of the first exit time $\tau (\gamma)$  introduced in Equation \ref{eq:FeynmanKacFormula} grows like $R^2$
(see Problem 2.25, page 253 of \cite{kara91a}).
In the first experiment of this section the potential will be computed along a horizontal line located $100$ m above the deposit. In the second experiment we compute the  potential along a vertical line intersecting the deposit. The situation is represented in Figure \ref{f1}.

We are now ready to present the results of the experiments obtained with two different implementations of the Monte Carlo method.
The first implementation is designed for the GPU: it is written in CUDA and tested on Neumann, a high-performance computer from Memorial University of Newfoundland, equipped with an NVIDIA Tesla K20Xm video card.
The second implementation is designed for multi-core parallel processing: it is written using the MPI protocol for C++, and tested on the compute cluster Torngat from Memorial University of Newfoundland, using 120 Intel Xeon X5660 cores.  
We report the performance and accuracy of the two versions of the algorithm. 
In the remainder of the paper we will denote the number of random walks simulated by the algorithm using $N$ and the time step of the algorithm described in Section 
\ref{sec:BackgroundMC} using $\Delta t$.
We remark that Neumann and Torngat will simulate over ranges of Monte Carlo simulations which are not exactly equal, but of similar magnitudes: this choice is due to the need of optimizing the usage of the computer architectures of the two machines.

\subsection{Experiment 1}

In the first numerical experiment we compute the values of the potential along the uniform grid of size of $10m$ on 
the horizontal line, where $y=50m$, $z=200m$, $x$ ranges from $-100m$ to $200m$ (with a grid spacing of $10m$).

\subsection{Comparison of performance: 1 GPU vs.\ 120 CPUs}

Table \ref{t3} and \ref{t4} report the execution time, at a given time step and number of Monte Carlo simulations, respectively for the CUDA and the MPI implementation. We compare the performance between one GPU and 119 parallel cores
(one of the 120 cores was used only to collect data from the other cores). From a comparison of the two tables we conclude that GPU wins in performance, as it is at least 2 times faster than using MPI on the traditional compute cluster. Moreover, since the MPI implementation runs in parallel among 119 processors, we conclude that the K20 Tesla GPU does the parallel job of about 200 compute cores.
Hence for the next experiments we chose to run simulations using the faster GPU cards.

\begin{table}[!ht]
	\caption{Experiment 1, GPU: execution time (measured in s) for $N$ Monte Carlo simulations and time step $\Delta t$.}
	\centering
	\begin{tabular}{ c| c c c c c c c c c}
		\hline	
		$\Delta t\downarrow$ / $ N\rightarrow$&  1024 & 2048 & 4096 & 5120 & 10240 & 20480 & 40960 & 51200\\
		\hline
		0.1 & 0.3 & 0.4 & 0.4 & 0.4 & 0.5 & 0.5 & 0.9 & 1.1\\
		0.05 & 0.8 & 0.9 & 1.0 & 1.0 & 1.1 & 1.2 & 2.0 & 2.5 \\
		0.01 & 4.8 & 5.3 & 5.8 & 5.9 & 6.4 & 7.1 & 11.5 & 14.3\\ 
		0.005 & 9.9 & 10.9 & 12.2 & 12.3 & 13.0 & 14.7 & 24.1 & 29.4\\
		0.001 & 52.8 & 55.6 & 59.6 & 61.1 & 65.7 & 75.8 & 125.6 & 155.7\\ 
		\hline
	\end{tabular}
	\label{t3}		
\end{table}

\begin{table}[!ht]
	\caption{Experiment 1, MPI: execution time (measured in $ 1s$) for $N$ Monte Carlo simulations and time step $\Delta t$, with 119 cores.}
	\centering
	\begin{tabular}{ c| c c c c c c c c c}
		\hline	
		$\Delta t\downarrow$ / $ N\rightarrow$&  
		1190 & 2380 & 3570  & 4760 & 11900 & 23800 & 35700  & 47600\\
		\hline
		0.1 & 0.1 & 0.2 & 0.3 & 0.4 & 0.8 & 1.5 & 2.3 & 2.9 \\
		0.05 & 0.3 & 0.4 & 0.6 & 0.8 & 1.7 & 3.2 & 4.7 & 6.3 \\ 
		0.01 & 1.3 & 2.2 & 3.2 & 4.0 & 9.2 & 17.6 & 26.1 & 34.2 \\
		0.005 & 2.6 & 4.6 & 6.7 & 8.2 & 19.1 & 36.3 & 53.2 & 70.9 \\ 
		0.001 & 13.6 & 23.7 & 34.3 & 43.0 & 98.7 & 189.4 & 278.1 & 365.0 \\
		\hline
	\end{tabular}	
	\label{t4}	
\end{table}

\subsection{Accuracy of the GPU}

\begin{figure}[!ht]
	\centering
	\begin{subfigure}[b]{0.5\textwidth}
		\centering
		\includegraphics[width=\textwidth]{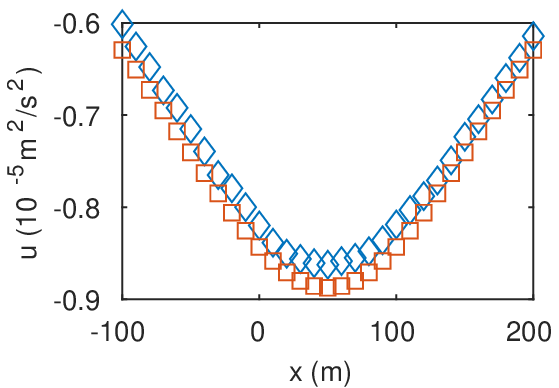}
	\end{subfigure}~
	\begin{subfigure}[b]{0.5\textwidth}
		\centering
		\includegraphics[width=\textwidth]{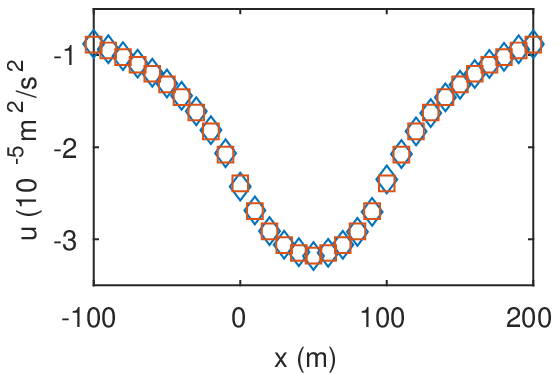}
	\end{subfigure}
	\caption{\textbf{Experiment 1:} The numerical results are obtained with time step $\Delta t=0.1$ and $N=51200$ Monte Carlo simulations \textit{(left)}. \textbf{Experiment 2:} The numerical results are obtained with time step $\Delta t=0.01$ and $N=51200$ Monte Carlo simulations \textit{(right)}.}
	\label{comparisons}\label{f2}
\end{figure}

In the left of Figure~\ref{comparisons} we depict the numerically computed values of the potential against the analytically computed values, along the $31$ grid points of the horizontal line located 100 m above the source ($z=200m$ with respect to our coordinate system).
Based on the previous considerations regarding the performance of the two different architectures, it is sufficient to report the results of the GPU simulation. A sufficiently small error was obtained at a time step $\Delta t=0.1$ and using $N=51200$ walks.  

The left panel of Figure \ref{comparisons} shows the translational constant error between the numerical and the analytical potentials, which is due to the incorrect choice of boundary conditions,
as explained in Section~\ref{subsec:boundary}. The RMS error is $1\%$,
and it was obtained with an execution time of $1$ s (as can be deduced from Table~\ref{t3}).

\subsection{Experiment 2}

In this experiment we find the numerical values of the potential along the vertical line  $y=50m$, $x=50m$, and $-100m\leqslant z \leqslant 200m$, on a uniform grid of size $\Delta z=10m$. 

We report that in order to obtain a RMS error comparable to the error obtained in Experiment 1, the time step had to be reduced to $\Delta t=0.01$.  
With $\Delta t=0.01$ we obtained an RMS error of $0.7\%$, while for $\Delta t=0.1 $ it was around 
$7\%$. The lower accuracy of this experiment is due to an increase in computational error for the potential at points close to or located at the boundary of the cube. The presence of a discontinuity surface of the source term affects the computations and moreover
shadows the systematic translational error due to the boundary conditions, as was observed in Experiment 1.  
The same issue will be discussed in Experiment 4 of Section \ref{sec:NumericalExperimentsMC1}, where we will show that the error
the acceleration computed as a finite difference derivative of the potential obtained with the Monte Carlo method, is difficult to reduce
for points close to the discontinuity surfaces of the source. 
We also remark that the translational error is also difficult to spot in Figure~\ref{f2} on the right panel since, 
as we move along a vertical line toward the source, the potential ranges over a larger scale with respect to which the translational error
is much smaller.

\section{Numerical experiments: physical results with multi GPU simulations}\label{sec:NumericalExperimentsMC1}

This set of experiments was simulated with a different implementation of the algorithm designed to run in parallel on a series of GPU cards. The new compute cluster at Memorial University consists of 8 nodes, each node having one NVIDIA Tesla K80 card. The Monte Carlo simulations were written in CUDA and executed on the GPUs. The GPUs were called independently across the nodes, hence in parallel, by the MPI based code. Each GPU computed several thousands of Monte Carlo simulations and the partial average over its sample space of random walks, according to the Feynman--Kac formula. Next, one node collected the partial results of each GPU and computed the final average as in \eqref{eq:Average} to produce the solution of the Poisson equation.

The possibility of using multiple GPUs opens to the perspective of dramatically reducing the time of computation. If provided with a large cluster of GPUs, a new level of parallelization can be introduced by letting different groups of GPUs compute the solution of the Poisson equation at different subsets of the grid. 

The code was written using floating point precision numbers, in order to minimize memory requirements and improve the performance. 
This choice however affects the precision: using floating point numbers, the relative error in the computation of the acceleration
could not be lowered by less than 1\%. However, this error being acceptable for a geophysical point of view, we decided to keep this precision for the next experiments.  We were able to achieve an error of 0.1\% using a collection of around 100 million of Monte Carlo simulations and a time step $\Delta t=0.1$, using double precision numbers rather than floating point.  

\subsection{Moving averages}\label{subsec:movav}

The Monte Carlo error in the calculation of the potential $u$ at the given set of points magnifies critically when the acceleration is computed using a finite difference approximation. The fluctuations on the results of the potential produce large deviations from the correct value of the derivative. In order to overcome such difficulty, we observe that the points $z_1,z_2,...z_n$, where the potential is obtained, are all located in a small region where the potential is  changing very slowly. Moreover, it is important to recall the values for the potential are obtained through fully independent simulations, that is, the values for the potential form a set of random continuous independent variables.

The sum of random variables allows one to devise a strategy to smooth the results obtained for the potential at the set of points $z_j, j=1,...,n$. In order to illustrate the strategy, consider as an example the sum of two independent normal random variables, that is, two variables whose values are normally distributed. The probability distribution of the sum of the two variables results from the convolution of the distributions. Namely, if $X$ and $Y$ are standard normal variables (whose density distributions are denoted by $f_X$ and $f_Y$) one has:

$$ f_X(x)=f_Y(y)=\frac{1}{\sqrt{2\pi}}e^{{-x^2}/2} \ .
$$

In this example, the convolution provides the distribution $f_Z$ for $Z=X+Y$, 
$$ f_Z(z)=f_X * f_Y(z)=\frac{1}{2\pi}\int_{-\infty}^{\infty}e^{{-(z-y)^2}/2}e^{{-y^2}/2} dy \ ,
$$
\noindent or, simply,
$$
f_Z(z)=\frac{1}{\sqrt{4\pi}}e^{{-z^2}/4} \ .
$$

The main observation is that there is a reduction of the variance, which supports the strategy of averaging values of the potential for a point, say, $z_i$, with its neighbors, aiming at a variance reduction and thus an improvement of the derivative to be computed.

Effectively, in order to reduce the error, we employed the method of moving averages. 
Specifically, the value of the potential $u$ at $z_k$ is then replaced by the average of the values of $u$ at a few neighboring points of $z_k$. We applied the method provided in \texttt{MATLAB} using the \texttt{smooth} function. As pointed out above,  the averaging process for the potential at a small collection of points helps reducing the variance. The outcomes from applying the method of moving averages were very successful as indicated by the results of the next experiment.

\subsection{Experiment 3}



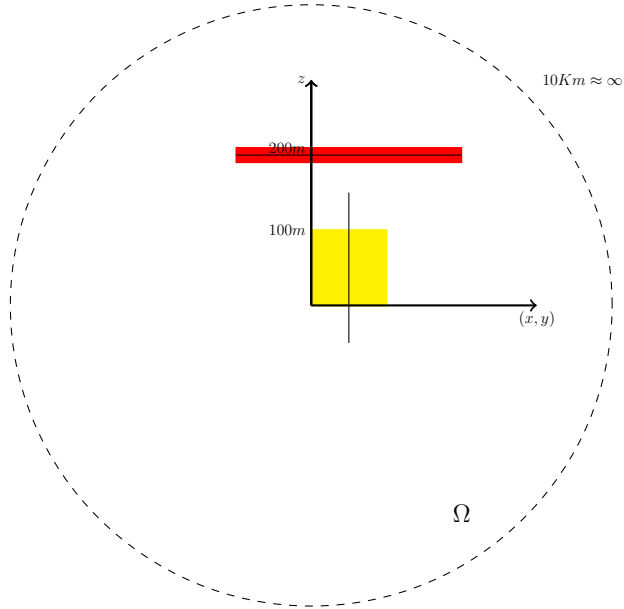
\begin{figure}
	\centering
	\begin{tikzpicture}
	\draw[dashed] (0,0) circle (4cm); 
	\filldraw[thick,yellow] (0,0)--(1,0)--(1,1)--(0,1)--cycle; 
	\filldraw[thin,red](-1,2.1)--(-1,1.9)--(2,1.9)--(2,2.1)--cycle; 
	\draw  (-1,2)--(2,2);
    \draw  (0.5,1.5)--(0.5,-0.5);
	\draw[thick,->] (0,0)--(0,3);
	\draw[thick, ->] (0,0)--(3,0) node[anchor=north,scale=0.5]{$(x,y)$};
	\draw[thick, ->] (0,0)--(0,3) node[anchor=east,scale=0.5]{$z$};
	\draw (0,1) node[anchor=east, scale=0.5]{$100m$};
	\draw (0,2.1) node[anchor=east, scale=0.5]{$200m$};	
	\draw (3,3) node[anchor=west,scale=0.5]{$10Km\approx\infty$};
	\draw (2,-3) node[anchor=south,scale=0.8]{$\Omega$};
	\end{tikzpicture}
	\caption{\textbf{Experiment 3:} Computation of the acceleration along the horizontal line. 
	\textbf{Experiment 4:} Computation of the acceleration along the vertical line.
	The yellow cube represents the mineral deposit.}
\end{figure}\label{f3}

The geometry of the problem is similar to the first set of experiments, with the main difference that here we chose a smaller vertical range and grid size for the rectangle where we compute the potential.  The potential is now computed on a 21 x 11 grid (total of 231 points). We reduced the  vertical spacing to of the grid points to $1m$, in order to find sufficiently accurate values for the vertical acceleration using a finite difference formula. 
The acceleration was computed, as in Experiment 1 along the horizontal line $y=50m$, $z=200m$, using the standard centered finite difference formula:  $u_z(x,z)=(u(x,z+\Delta z)-u(x,z-\Delta z))/(2\Delta z)$,
where $\Delta z=1m$.

We remark that there exists a closed analytic formula for the vertical acceleration
\cite{blakely1996potential}. However, it is rather complicated and at points above the boundary of the source
assumes an indeterminate form. Hence we decided to use a finite difference computation. 
Figure \ref{f3.5} compare the calculations for the acceleration along the horizontal line computed
in the two ways, showing a very good agreement. 

\begin{figure}
\centering
\includegraphics[width=0.5\textwidth]{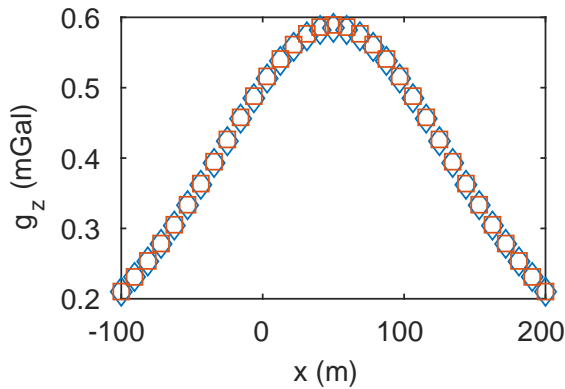}
\caption{\textbf{Experiment 3:} computation of the vertical acceleration $g_z$ using the analytic formula (diamonds)
	and using a finite difference approximation (squares).}
\end{figure}\label{f3.5}

\begin{figure}[!ht]
	\centering
	\begin{subfigure}[b]{0.5\textwidth}
		\centering
		\includegraphics[width=\textwidth]{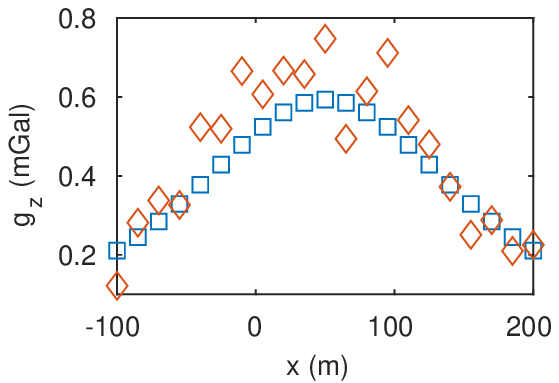}
	\end{subfigure}~
	\begin{subfigure}[b]{0.5\textwidth}
		\centering
		\includegraphics[width=\textwidth]{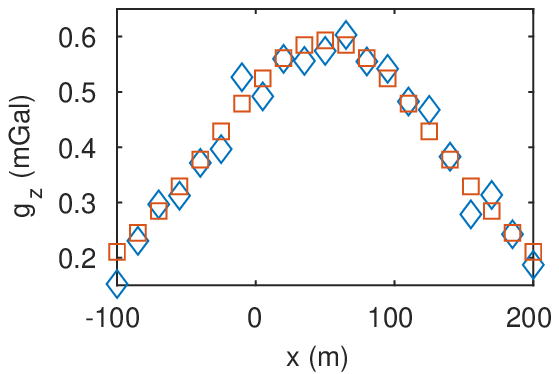}
	\end{subfigure}
	\caption{\textbf{Experiment 3:} computation of the acceleration. \textit{Left:} without smoothing (RMS error=0.21); \textit{Right:} with smoothing (RMS error=0.06). The potential is computed with 8 GPUs running in parallel(time step $ \Delta t$ =0.1, N=65536 Monte Carlo simulations, execution time= 7s)}
		\label{f4}
\end{figure}

\begin{figure}[!ht]
	\centering
	\begin{subfigure}[b]{0.5\textwidth}
		\centering
		\includegraphics[width=\textwidth]{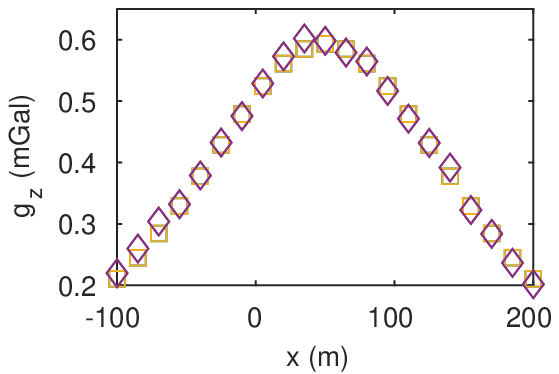}
	\end{subfigure}~
	\begin{subfigure}[b]{0.5\textwidth}
		\centering
		\includegraphics[width=\textwidth]{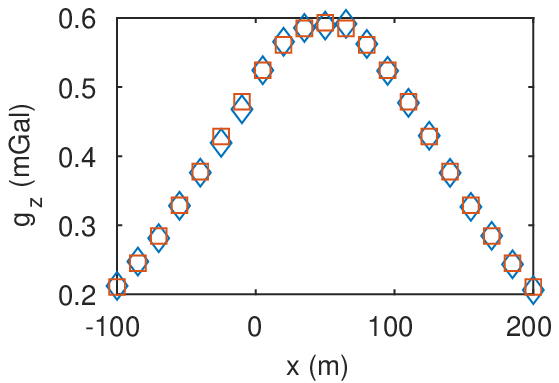}
	\end{subfigure}
	\caption{\textbf{Experiment 3:} computation of the acceleration. \textit{Left:} without smoothing (RMS error=0.02). \textit{Right:} with smoothing (RMS error=0.009). The potential is computed with 8 GPUs running in parallel
		(time step $\Delta t$= 0.1, N = 4194304 Monte Carlo simulations, execution time= 31s}
		\label{f5}
\end{figure}

Figures \ref{f4} and \ref{f5} show the improvement of the calculation of the acceleration obtained by applying the method of moving averages. The method is more advantageous when the number of random walks performed is relatively small, as shown in 
\ref{f5}: here the smoothing procedure produces a drastic improvement in the computation of the acceleration. 
With a larger collection of walks, as shown in the graphs in Figure \ref{f6}, the improvement is less significant. 
This makes the method of moving average a powerful tool since it allows us to reduce the number of Monte Carlo simulations
necessary to obtain an acceptable accuracy.

\subsection{Comparison with finite element methods}\label{sub:FE-comparison}

\begin{figure}[!ht]
	\centering
	\includegraphics[width=0.5\textwidth]{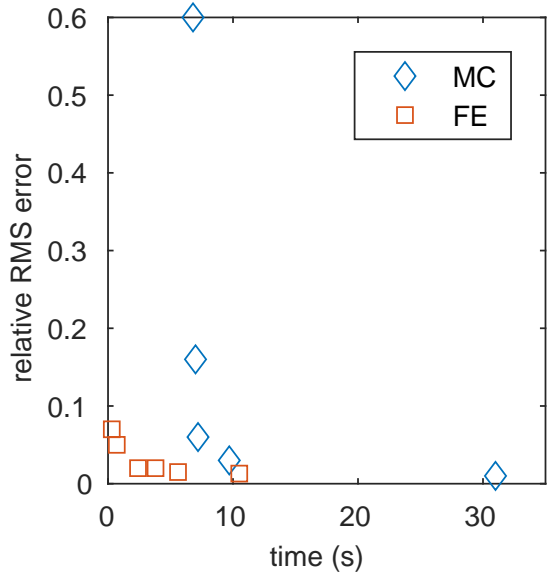}
	\caption{Comparison: Finite Elements (FE) vs.\ Monte Carlo (MC).}	
	\label{f6}
\end{figure}

Figure~\ref{f6} presents a comparison of the performances of the Monte Carlo algorithm with a 
deterministic algorithm, based on finite elements~\cite{jaha13a}. 
The plot shows that the two methods tie in performance at an execution time of the order of 10 seconds. However, our parallel implementation of the Monte Carlo method was tested on a machine endowed with 8 GPUs only. More powerful machines are available in scientific computing labs, and those can be used to not only tie, but beat the performance of a deterministic method. 
We infer the following estimate. We computed the solution of the Poisson equation at 126 points, using an algorithm where each GPU card computed the same number of Monte Carlo simulations at all points of the grid. With 80 GPUs, the job could be subdivided at the level of points as well: we could assign to a group of 10 GPUs the computation of the solution for a subset of the grid that consists of about 12 points (126/10): that will reduce the execution time by a factor of 10, hence beating the finite element method, provided that we are satisfied with an error of $1\%$, acceptable in the contest of mine exploration. Compute machines with about 100 GPUs are already available in some computational center, which lends credibility to our estimate using 80 GPUs. For example, very recently the College of Engineering of Carnegie Mellon University installed a GPU cluster machine endowed with 112 GPUs\footnote{\url{https://www.cmu.edu/me/news/archive/2017/viswanathan-launches-gpu-cluster.html}}. 

 We were able to achieve an error of 0.1\% using a collection of around 100 milion of Monte Carlo simulations and a time step
$\Delta t=0.1$, using double precision numbers rather than floating point.

\subsection{Experiment 4: Discontinuity of the source}
As we mentioned in Section~\ref{subsec:boundary}, the stochastic solution of the Poisson equation requires
the source term to be at least continuous. However, in all our experiments the source has jump discontinuities along the boundary of the cube. We show here that the standard Feynman--Kac formula still holds approximately for this case, and gives approximation errors that are well within the geophysically required tolerance. The investigation of an adjustment of this formula to the case of discontinuous source terms will be considered elsewhere.

\begin{figure}[!ht]
	\centering
	\begin{subfigure}[b]{0.5\textwidth}
		\centering
		\includegraphics[width=\textwidth]{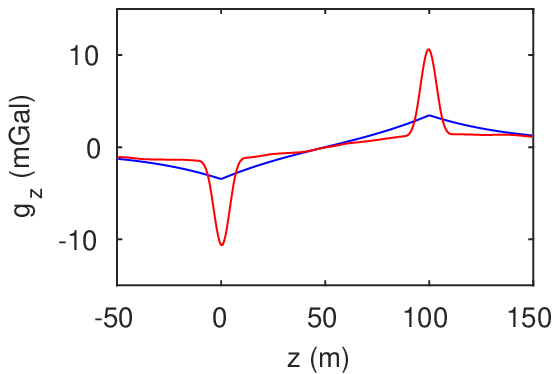}
	\end{subfigure}~
	\begin{subfigure}[b]{0.5\textwidth}
		\centering
		\includegraphics[width=\textwidth]{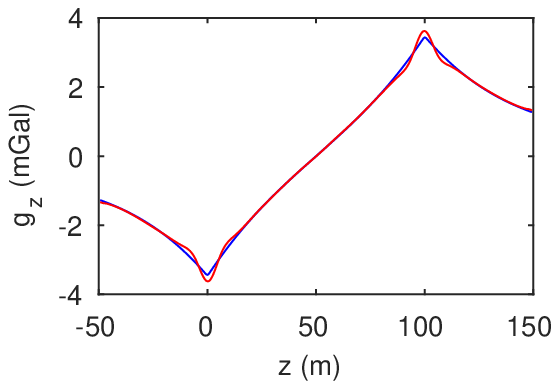}
	\end{subfigure}
	\caption{\textbf{Experiment 4:} Computation of the acceleration along a vertical line crossing the mineral deposit. The red curve represents
		the acceleration computed from the Monte Carlo simulation, the blue curve the acceleration computed as finite difference derivative
		of the exact potential. The left panel is obtained with a time step $\Delta t =0.1 $, $N=65536$ Monte Carlo simulations; the right panel is obtained with $\Delta t =0.01$ and $N=4194304$ Monte Carlo simulations.}
\label{f8}
\end{figure}

Here, we analyze the effect on the discontinuity by computing the acceleration along a vertical line that crosses the cube. The graphs plotted in Figure~\ref{f8} depict clearly that the worst discrepancy between the numerical and the analytical values of the acceleration occur at the boundary of the cube. We note however that we observed that even for smooth sources the high discrepancy occurs at inflection points of the potential.   

\section{Conclusions}\label{sec:ConclusionsMC}

We showed that Monte Carlo based algorithms can be a powerful alternative to deterministic methods for solving model problems in geophysics. In this paper we consider only simple cases, while in a future paper we will analyze more complex situations, extending our research to the solution of the two-dimensional Maxwell's equations.
Table \ref{table:conclusion} reports some comparison between Monte Carlo methods and Finite Elements methods (or, for that matter, other grid based deterministic methods). We conclude by highlighting some key results of our numerical experiments.

\begin{table}[!ht]
	\caption{Comparison of stochastic methods with finite element based methods.}
	\centering
	\begin{tabular}{  |p{4cm}|p{5cm}|p{5cm}| }
		\hline
		& {\bf Stochastic methods} & { \bf Finite Elements}\\
		\hline
		{\bf Speed} & proportional to the number of  points at which the solution is computed & depends on the grid size\\
		\hline
		{\bf Memory} & inexpensive (does not require a grid) & expensive (requires a grid)\\
		\hline
		{\bf Speed vs.\ accuracy} & fast performance for a $1\%$ error & fast performance for a $1\%$ error
		\\
		\hline	
		{\bf Implementation} & very easy, suitable for massive parallelization & not as easily parallelizable\\
		\hline	
	\end{tabular}	
	\label{table:conclusion}	
\end{table}

{ \bf Speed and Accuracy.}
For typical geophysical accuracy requirements of about 1\% error, the stochastic method performs exceptionally well with a limited number of Monte Carlo simulations that do not require the existence of a particularly large compute cluster. Higher accuracy is achievable as well using sufficiently many Monte Carlo simulations with fine enough time steps, which are easily obtained with a moderately large high performance compute cluster.

\smallskip

{\bf Implementation.}
The design of an algorithm that can run parallel Monte Carlo simulations is rather simple, as shown in Section~\ref{sec:BackgroundMC}. Since the stochastic method is embarrassingly parallel, not much effort has to be spent on writing a reliable, GPU enabled code to implement this method, making it readily implementable for commercial applications.

\subsection*{Acknowledgements}

This research was undertaken, in part, thanks to funding from the Canada Research Chairs program, the NSERC Discovery Grant program and the IgniteR{\&}D program of the Research and Development Corporation of Newfoundland and Labrador (RDC).
The authors would like to thank the following contributers, for their input, technical support
and expertise: Lawrence Greening, Hormoz Jahandari, J.-C. Loredo-Osti, Paul Sherren, Lukas Spies, Oliver Stueker, Lara Zabel.

\footnotesize\setlength{\itemsep}{0ex}
\bibliographystyle{siam}

\begin{thebibliography}{10}
	
	\bibitem{aceb05a}
	{\sc J.~A. Acebr{\'o}n, M.~P. Busico, P.~Lanucara, and R.~Spigler}, {\em Domain
		decomposition solution of elliptic boundary-value problems via monte carlo
		and quasi-monte carlo methods}, SIAM J. Sci. Comput., 27 (2005),
	pp.~440--457.
	
	\bibitem{aceb10a}
	{\sc J.~A. Acebr{\'o}n, {\'A}.~Rodr{\'\i}guez-Rozas, and R.~Spigler}, {\em
		Efficient parallel solution of nonlinear parabolic partial differential
		equations by a probabilistic domain decomposition}, J. Sci. Comput., 43
	(2010), pp.~135--157.
	
	\bibitem{aceb07a}
	{\sc J.~A. Acebr{\'o}n and R.~Spigler}, {\em A new probabilistic approach to
		the domain decomposition method}, in Domain Decomposition Methods in Science
	and Engineering XVI, Springer, 2007, pp.~473--480.
	
	\bibitem{bihl16b}
	{\sc A.~Bihlo, C.~G. Farquharson, R.~D. Haynes, and J.~C. Loredo-Osti}, {\em
		Probabilistic domain decomposition for the solution of the two-dimensional
		magnetotelluric problem}, Computat. Geosci.,  (2016), pp.~1--13.
	
	\bibitem{bihl14a}
	{\sc A.~Bihlo and R.~D. Haynes}, {\em Parallel stochastic methods for pde based
		grid generation}, Comput. Math. Appl., 68 (2014), pp.~804--820.
	
	\bibitem{bihl15b}
	{\sc A.~Bihlo, R.~D. Haynes, and E.~J. Walsh}, {\em Stochastic domain
		decomposition for time dependent adaptive mesh generation}, J. Math. Study,
	48 (2015), pp.~106--124.
	
	\bibitem{blakely1996potential}
	{\sc R.~J. Blakely}, {\em Potential theory in gravity and magnetic
		applications}, Cambridge University Press, 1996.
	
	\bibitem{farquharson2009three}
	{\sc C.~Farquharson and C.~Mosher}, {\em Three-dimensional modelling of gravity
		data using finite differences}, Journal of Applied Geophysics, 68 (2009),
	pp.~417--422.
	
	\bibitem{gobe00a}
	{\sc E.~Gobet}, {\em Weak approximation of killed diffusion using euler
		schemes}, Stochastic Process. Appl., 87 (2000), pp.~167--197.
	
	\bibitem{jaha13a}
	{\sc H.~Jahandari and C.~G. Farquharson}, {\em Forward modeling of gravity data
		using finite-volume and finite-element methods on unstructured grids},
	Geophysics, 78 (2013), pp.~G69--G80.
	
	\bibitem{kara91a}
	{\sc I.~Karatzas and S.~E. Shreve}, {\em Brownian motion and stochastic
		calculus}, vol.~113 of Graduate Texts in Mathematics, Springer, New York,
	1991.
	
	\bibitem{mils04a}
	{\sc G.~Milstein and M.~Tretyakov}, {\em Stochastic numerics for mathematical
		physics}, Springer, Berlin Heidelberg, 2004.
	
	\bibitem{murray2001best}
	{\sc A.~S. Murray and R.~M. Tracey}, {\em Best practices in gravity surveying},
	GeoscienceAustralia, 3 (2001).
	
	\bibitem{pres07Ay}
	{\sc W.~H. Press, S.~A. Teukolsky, W.~T. Vetterling, and B.~P. Flannery}, {\em
		Numerical recipes 3rd edition: {T}he art of scientific computing}, Cambridge
	University Press, Cambridge, UK, 2007.
	
	\bibitem{waldvogel1976newtonian}
	{\sc J.~Waldvogel}, {\em The newtonian potential of a homogeneous cube},
	Zeitschrift f{\"u}r Angewandte Mathematik und Physik (ZAMP), 27 (1976),
	pp.~867--871.
	
\end{thebibliography}

\end{document}